\documentclass{amsart}
\usepackage[T1]{fontenc}
\usepackage{amsmath,amssymb,amsthm,mathtools}
\usepackage{a4wide}
\usepackage[colorlinks=true,linkcolor=blue,citecolor=blue,urlcolor=blue]{hyperref}
\newtheorem{theorem}{Theorem}[section]
\newtheorem{lemma}[theorem]{Lemma}
\newtheorem{proposition}[theorem]{Proposition}
\newtheorem{corollary}[theorem]{Corollary}
\numberwithin{equation}{section}
\newcommand{\Zc}{Z_{\mathrm c}}
\newcommand{\core}{\mathrm c}
\newcommand{\U}{\mathcal U}
\newcommand{\E}{\mathbb E}
\title[Banach isomorphism between hyperfinite type III \(L_p\)-spaces]
{Banach isomorphism between separable \(L_p\)-spaces, \(1\le p\ne2<\infty\),
associated with hyperfinite type III von Neumann algebras}
\author[C. Gao]{Chenhao Gao}
\address{Institute for Advanced Study in Mathematics of HIT, Harbin Institute
  of Technology, Harbin, 150001, China}
\author[J. Huang]{Jinghao Huang}
\address{Institute for Advanced Study in Mathematics of HIT,
Harbin Institute of Technology, Harbin, 150001, China}
\thanks{J. Huang was supported by the National Natural Science Foundation
of China under grant nos.~12301160, 12471134 and 12671159.}
\author[Y. Zhu]{Yunpeng Zhu}
\address{Institute for Advanced Study in Mathematics of HIT, Harbin Institute
  of Technology, Harbin 150001, China}
\date{}
\begin{document}
\begin{abstract}
For each fixed \(p\in[1,\infty)\setminus\{2\}\), the number of Banach
isomorphism classes of separable \(L_p\)-spaces associated with
hyperfinite type III von Neumann algebras is \(2^{\aleph_0}\).
\end{abstract}
\maketitle

\section{Introduction}

For each fixed \(1\le p<\infty\), \(p\ne2\), the classical separable
infinite-dimensional \(L_p(\mu)\)-spaces fall into exactly two Banach
isomorphism classes, represented by \(\ell_p\) and \(L_p(0,1)\)
\cite[Section~4.9.3.3, p.~136]{Pietsch2007}.
This raises the corresponding question for the spaces \(L_p(M)\)
associated with hyperfinite von Neumann algebras \(M\).
Sukochev~\cite{Sukochev2000} studied this question for algebras of
type~I. Haagerup, Rosenthal and Sukochev~\cite[Theorem~5.1]{HRS}
gave a complete classification in the semifinite case. They proved
that the separable infinite-dimensional \(L_p\)-spaces associated
with hyperfinite semifinite von Neumann algebras fall into exactly
thirteen Banach isomorphism classes.

Huang and Sukochev~\cite[Theorem~1.1]{HuangSukochev2026} showed that
none of these thirteen classes contains \(L_p(M)\) for a nonzero
type~III von Neumann algebra \(M\) with separable predual.
More generally, for any semifinite von Neumann algebra \(N\), they
proved that \(M_*\) does not embed isomorphically into \(N_*\) and
that, for \(1<p<\infty\), \(p\ne2\), \(L_p(M)\) is not isomorphic
to a complemented subspace of \(L_p(N)\).

For hyperfinite type~III factors with separable preduals,
Haagerup, Rosenthal and Sukochev~\cite[Theorem~6.5, p.~65]{HRS} proved
that \(L_p(R_\lambda)\) is completely isomorphic to \(L_p(R_\infty)\)
for every \(0<\lambda<1\), where \(R_\lambda\) and \(R_\infty\)
denote the hyperfinite factors of types \(\mathrm{III}_\lambda\)
and \(\mathrm{III}_1\), respectively. They also exhibited uncountably
many pairwise nonisomorphic hyperfinite \(\mathrm{III}_0\) factors
whose \(L_p\)-spaces all belong to this same Banach isomorphism class.
At \(p=1\), Haagerup and Musat~\cite[Theorems~4.1 and~4.2]{HM}
constructed a family \((M^{(t)})_{0\le t<2}\) of hyperfinite
\(\mathrm{III}_0\) factors such that, for \(s\ne t\), there is no
completely bounded isomorphism between the preduals
\((M^{(s)})_*\) and \((M^{(t)})_*\).

In this paper we consider isomorphisms of the underlying Banach spaces.
Theorem~\ref{thm:criterion} characterizes Banach isomorphism of
\(L_p(M)\) and \(L_p(N)\) in terms of the flows of weights of \(M\)
and \(N\). Applying this criterion to the Haagerup--Musat family,
we obtain \(L_p(M^{(s)})\not\simeq L_p(M^{(t)})\) whenever
\(s\ne t\). Thus, for each fixed \(1\le p<\infty\), \(p\ne2\),
exactly \(2^{\aleph_0}\) Banach isomorphism classes occur in the
hyperfinite type~III case, and this cardinality is already attained
among factors of type \(\mathrm{III}_0\)
(Corollary~\ref{cor:continuum-classes}).

For a von Neumann algebra \(M\), let \(L_p(M)\) denote its Haagerup \(L_p\)-space.
For Banach spaces \(X,Y\), write \(X\hookrightarrow Y\) if there is a
bounded linear isomorphic embedding of \(X\) into \(Y\), and
\(X\simeq Y\) if \(X\) and \(Y\) are Banach isomorphic.
For von Neumann algebras, \(\cong\) denotes a normal
\(*\)-isomorphism. Write \(\U(M)\) for the unitary group of \(M\).
For a positive normal functional \(\varphi\in M_*^+\), let
\(s(\varphi)\in M\) be its support projection and let
\(\varphi^{1/p}=D_\varphi^{1/p}\in L_p(M)_+\), where \(D_\varphi\)
is its Haagerup density. For a faithful normal state \(\varphi\),
write \(\sigma^\varphi\) for its modular automorphism group and
\[
 M_\varphi=\{a\in M:\varphi(ax)=\varphi(xa)\text{ for all }x\in M\}
\]
for its centralizer.

Let \(\core(M)\) be the continuous core of \(M\), and set
\(\Zc(M)=Z(\core(M))\). The restriction of the dual action to
\(\Zc(M)\) is denoted by \(\theta^M\). The pair
\(\mathcal F(M)=(\Zc(M),\theta^M)\) is the flow of weights of \(M\).
We write \(\mathcal F(N)\to\mathcal F(M)\) if there is a normal unital
positive linear map
\begin{equation}\label{eq:flow-map-definition}
 K:\Zc(N)\longrightarrow\Zc(M),\qquad
 K\circ\theta^N_t=\theta^M_t\circ K\quad(t\in\mathbb R).
\end{equation}
A map satisfying \eqref{eq:flow-map-definition} is called an equivariant
Markov map.

\begin{theorem}\label{thm:criterion}
Let \(M,N\) be nonzero hyperfinite von Neumann algebras of pure type III
with separable preduals.
\begin{enumerate}
\item For \(1\le p<2\), we have \(L_p(M)\hookrightarrow L_p(N)\) if and
only if \(\mathcal F(N)\to\mathcal F(M)\). If
\(\mathcal F(N)\to\mathcal F(M)\), then \(L_p(M)\) is isometric to a
\(1\)-complemented subspace of \(L_p(N)\) for every \(1\le p<\infty\).
\item For each \(1\le p<\infty,\ p\ne2\), the spaces
\(L_p(M)\) and \(L_p(N)\) are Banach isomorphic if and only if
\(\mathcal F(N)\to\mathcal F(M)\) and
\(\mathcal F(M)\to\mathcal F(N)\).
\end{enumerate}
\end{theorem}

\section{Recovering a map of flows from a Banach embedding}

\subsection{Maps of flows associated with tensor products, ultrapowers and expectations}
\label{subsec:ambient-transfer}
Throughout Section~2, let \(M,N\) be as in
Theorem~\ref{thm:criterion}. All ultrapowers are taken along free
ultrafilters on countable sets. We first describe the maps of flows
used to construct the target algebra of an isometric embedding.
Let \(B\) be a properly infinite von Neumann algebra with separable
predual.
\begin{enumerate}
\item For \(B'=L_\infty(0,1)\bar\otimes B\), tensoring with the tracial
abelian algebra gives
\[
 \core(B')\cong L_\infty(0,1)\bar\otimes\core(B).
\]
Under this identification, the map
\[
 \Zc(B)\longrightarrow\Zc(B'),\qquad z\longmapsto1\otimes z
\]
is a normal unital equivariant \(*\)-homomorphism.
\item Let \(\psi\) be a faithful normal state on \(B\), and let
\(B^\omega\) be its Ocneanu ultrapower, with ultrapower state
\(\psi^\omega\). Choose a von Neumann subalgebra \(B'\subset B^\omega\)
with separable predual that contains the diagonal copy of \(B\) and
is invariant under \(\sigma^{\psi^\omega}\).
Tomatsu's crossed-product embedding \cite[Theorem~3.6]{Tomatsu}, applied
to the modular action, realizes \(\core(B^\omega)\) inside \(\core(B)^\omega\)
and preserves the diagonal copy of \(\core(B)\). Here the modularly
equicontinuous part of \(B^\omega\) is all of \(B^\omega\)
\cite{AH,Tomatsu}. For \(z\in\Zc(B)\), its constant representative
belongs to \(\core(B^\omega)\) and commutes with \(\core(B)^\omega\).
Thus, the diagonal inclusion gives a normal unital equivariant
\(*\)-homomorphism \(\Zc(B)\to\Zc(B^\omega)\).
Takesaki's theorem gives a faithful normal
\(\psi^\omega\)-preserving conditional expectation
\(B^\omega\to B'\). Restricting the induced core expectation to the
center, as described in item~3, gives
\(\mathcal F(B)\to\mathcal F(B')\). The algebra \(B'\) is properly
infinite because it contains the diagonal copy of \(B\).
\item Let \(P\subset W\) be an inclusion of \(\sigma\)-finite von
Neumann algebras with faithful normal conditional expectation
\(E:W\to P\). Choose a faithful normal state \(\gamma\) on \(P\), and
use \(\gamma\circ E\) to form the core of \(W\).
The induced conditional expectation satisfies
\[
 \widetilde E:\core(W)\longrightarrow\core(P),\qquad
 \widetilde E(a\lambda_t)=E(a)\lambda_t
 \quad(a\in W,\ t\in\mathbb R),
\]
where \(\lambda_t\) denotes the canonical implementing unitary in the
respective cores. The map \(\widetilde E\) intertwines the dual actions.
Its \(\core(P)\)-bimodularity gives
\(\widetilde E(\Zc(W))\subset\Zc(P)\), and hence
\(\mathcal F(W)\to\mathcal F(P)\)
\cite[proof of Theorem~3.5]{HM}.

For a nonzero projection \(f\in W\), choose a faithful normal state
on \(W\) whose centralizer contains \(f\). Then
\(\core(fWf)\cong f\core(W)f\), and the normal unital equivariant map
\[
 \Zc(W)\longrightarrow Z(f\core(W)f),\qquad z\longmapsto fz
\]
gives \(\mathcal F(W)\to\mathcal F(fWf)\).
\end{enumerate}

We also need to realize certain subspaces of the Banach ultrapower
\(L_p(B)^\omega\) in the \(L_p\)-space of an algebra \(B'\) as in
item~2. Let \(1\le p<\infty\). Following Raynaud and Xu, a bounded
set \(\mathcal K\subset L_p(B)\) is \(p\)-equiintegrable if
\(\sup_{y\in\mathcal K}\|f_\alpha y f_\alpha\|_{L_p(B)}\to0\)
for every net of projections \(f_\alpha\in B\) converging ultraweakly
to zero \cite{RX}. For a bounded family, \(p\)-equiintegrability
means \(p\)-equiintegrability of its range.

Let \(G\) be the Groh--Raynaud ultrapower in Raynaud's identification
\(L_p(B)^\omega=L_p(G)\) \cite{Raynaud}, and let \(e\in G\) be the
support of the constant class of \(\psi\in B_*^+\). Then
\(eGe\cong B^\omega\) \cite{AH,RX}. Write \(z_G(e)\) for the
central support of \(e\) in \(G\). If \(y\in L_p(G)\) admits a
\(p\)-equiintegrable representative, then
\((1-e)y(1-e)=0\) by \cite[Theorem~4.6]{RX}, and therefore
\[
                 y=z_G(e)y z_G(e).
\]
Let \(E\subset L_p(G)\) be a separable subspace whose elements
satisfy this support identity. Choose a dense sequence
\((y_m)_{m\ge1}\) in \(E\), and let \(r\) be the join of \(e\)
and the left and right support projections of all \(y_m\).
The projection \(r\) is \(\sigma\)-finite and \(r\le z_G(e)\).
Moreover, \(r\precsim e\), where \(\precsim\) denotes
Murray--von Neumann subequivalence. To see this, take a maximal
orthogonal family of subprojections of \(r\), each subequivalent to
\(e\). The family is countable because \(r\) is \(\sigma\)-finite.
A nonzero remainder \(r_0\le z_G(e)\) would satisfy \(r_0Ge\ne0\).
The polar decomposition of a nonzero element of \(r_0Ge\) would
give a further nonzero subprojection of \(r_0\) subequivalent to
\(e\), contrary to maximality. Since \(eGe\) is properly infinite,
the countable orthogonal sum is subequivalent to \(e\).
Choose a partial isometry \(v\in G\) with \(v^*v=r\) and \(vv^*\le e\).
We obtain the linear isometry
\begin{equation}\label{eq:ambient-support-transfer}
 E\longrightarrow L_p(eGe)=L_p(B^\omega),\qquad y\longmapsto vyv^*.
\end{equation}

To choose \(B'\subset B^\omega\) with \(vEv^*\subset L_p(B')\),
approximate each \(vy_mv^*\) in \(L_p(B^\omega)\) by a sequence
of elements \(a_{m,k}(\psi^\omega)^{1/p}\), with
\(a_{m,k}\in B^\omega\). Generate \(B'\) by the diagonal copy of
\(B\), all \(a_{m,k}\), and their modular translates at rational
times. Strong continuity of \(\sigma^{\psi^\omega}\) makes \(B'\)
invariant under the modular group. Takesaki's theorem gives a
faithful normal state-preserving conditional expectation
\(B^\omega\to B'\), and \(vEv^*\subset L_p(B')\) by the chosen
approximations. The same construction permits the inclusion of any
prescribed countable family of elements of \(B^\omega\) in \(B'\).

\subsection{Averaging an arbitrary embedding}\label{subsec:averaging}
Fix \(1\le p<2\). Suppose that
\(T:L_p(M)\hookrightarrow L_p(N)\), and choose constants
\(0<c\le C<\infty\) such that
\[
 c\|x\|_{L_p(M)}\le\|Tx\|_{L_p(N)}
       \le C\|x\|_{L_p(M)}\qquad(x\in L_p(M)).
\]
Let \((M_n)_{n\ge1}\) be an increasing sequence of finite-dimensional
unital von Neumann subalgebras whose union is strongly dense in \(M\).
For each \(n\), choose a measurable map
\(u_n:(0,1)\to\U(M_n)\) that sends Lebesgue measure to normalized
Haar measure on \(\U(M_n)\). Put \(D=L_\infty(0,1)\), and define
\begin{equation}\label{eq:average-map}
 \begin{gathered}
 F_n:L_p(M)\longrightarrow L_p(D\bar\otimes N)
                =L_p((0,1);L_p(N)),\\
 (F_nx)(t)=T(u_n(t)x)\qquad(x\in L_p(M),\ t\in(0,1)).
 \end{gathered}
\end{equation}
For every \(x\in L_p(M)\), integration gives
\begin{equation}\label{eq:average-embedding-bounds}
 c\|x\|_{L_p(M)}\le\|F_nx\|_{L_p(D\bar\otimes N)}
                 \le C\|x\|_{L_p(M)}.
\end{equation}

We prove that \(\{F_nx:n\ge1\}\) is \(p\)-equiintegrable in
\(L_p(D\bar\otimes N)\) for each fixed \(x\in L_p(M)\).
Choose \(r>0\) with \(1/r=1/p-1/2\). The polar decomposition of
\(x\) gives a factorization \(x=ab\), with \(a\in L_2(M)\),
\(b\in L_r(M)\), and
\(\|a\|_{L_2(M)}\|b\|_{L_r(M)}=\|x\|_{L_p(M)}\).
Let \(\varepsilon_1,\ldots,\varepsilon_k\) be independent random
signs, each with probabilities \(1/2,1/2\). Thus, \(\E\) denotes
\(2^{-k}\sum_{\varepsilon\in\{-1,1\}^k}\).
For any \(n_1,\ldots,n_k\ge1\), H\"older's inequality and
the Hilbert-space identity for the second sign moment give
\begin{equation}\label{eq:average-rademacher-bound}
\begin{aligned}
 \E\left\|\sum_{\ell=1}^k\varepsilon_\ell F_{n_\ell}x
                      \right\|_{L_p(D\bar\otimes N)}^p
 &\le C^p\|b\|_{L_r(M)}^p
       \int_0^1\E\left\|\sum_{\ell=1}^k
            \varepsilon_\ell u_{n_\ell}(t)a\right\|_{L_2(M)}^p\,dt\\
 &\le C^p\|b\|_{L_r(M)}^p
       \int_0^1\left(\sum_{\ell=1}^k
                    \|u_{n_\ell}(t)a\|_{L_2(M)}^2\right)^{p/2}\,dt\\
 &=C^p k^{p/2}\|x\|_{L_p(M)}^p.
\end{aligned}
\end{equation}
The second inequality uses \(p<2\).
The same factorization, with contractions \(a_\ell\in M\) in place
of \(u_{n_\ell}(t)\), proves the first estimate below. Factoring
\(x\) as an element of \(L_r(M)L_2(M)\) proves the second estimate.
For every \(x\in L_p(M)\) and \(\|a_\ell\|_\infty\le1\), we have
\begin{equation}\label{eq:one-sided-rademacher-bound}
 \begin{aligned}
 \left(\E\left\|\sum_{\ell=1}^k
             \varepsilon_\ell T(a_\ell x)\right\|_{L_p(N)}^p\right)^{1/p}
     &\le C\sqrt{k}\,\|x\|_{L_p(M)},\\
 \left(\E\left\|\sum_{\ell=1}^k
             \varepsilon_\ell T(xa_\ell)\right\|_{L_p(N)}^p\right)^{1/p}
     &\le C\sqrt{k}\,\|x\|_{L_p(M)}.
 \end{aligned}
\end{equation}

To deduce equiintegrability from these estimates, let
\(\mathcal K\subset L_p(W)\) be a bounded family, where \(W\) is
any von Neumann algebra. If \(\mathcal K\) fails
\(p\)-equiintegrability, \cite[Corollary~4.9]{RX} gives pairwise
orthogonal projections \(f_\ell\in W\), elements
\(y_\ell\in\mathcal K\), and \(\delta>0\) such that
\(\|f_\ell y_\ell f_\ell\|_{L_p(W)}\ge\delta\).
For \(k\ge1\), let
\[
 P_k(y)=\sum_{\ell=1}^k f_\ell y f_\ell,\qquad
 w_\varepsilon=\sum_{\ell=1}^k\varepsilon_\ell f_\ell
                     +1-\sum_{\ell=1}^k f_\ell.
\]
The map \(P_k:L_p(W)\to L_p(W)\) is contractive, and
\(w_\varepsilon\in\U(W)\). Consequently,
\begin{equation}\label{eq:equiintegrability-sign-obstruction}
\begin{aligned}
 \delta k^{1/p}
 &\le\left\|\sum_{\ell=1}^k f_\ell y_\ell f_\ell\right\|_{L_p(W)}\\
 &=\left\|\E\left[w_\varepsilon P_k
              \left(\sum_{\ell=1}^k\varepsilon_\ell y_\ell\right)
          \right]\right\|_{L_p(W)}
 \le\E\left\|\sum_{\ell=1}^k\varepsilon_\ell y_\ell\right\|_{L_p(W)}.
\end{aligned}
\end{equation}
Since \(1/p>1/2\), \eqref{eq:equiintegrability-sign-obstruction}
contradicts \eqref{eq:average-rademacher-bound} when
\(\mathcal K=\{F_nx:n\ge1\}\). The same argument applied to
\eqref{eq:one-sided-rademacher-bound} proves equiintegrability of
\(\{T(ax):\|a\|_\infty\le1\}\) and
\(\{T(xa):\|a\|_\infty\le1\}\).

Fix a free ultrafilter \(\omega\) on \(\mathbb N\). Apply
\eqref{eq:ambient-support-transfer}, with \(B=D\bar\otimes N\),
to the range of \(x\mapsto(F_nx)^\omega\). We obtain a properly
infinite von Neumann algebra \(B_0\) with separable predual and a
bounded linear map \(J_0:L_p(M)\to L_p(B_0)\) such that
\[
 \|J_0x\|_{L_p(B_0)}
     =\lim_{n\to\omega}\|F_nx\|_{L_p(D\bar\otimes N)}
 \quad(x\in L_p(M)),\qquad
 \mathcal F(N)\to\mathcal F(B_0).
\]
In particular, \(J_0\) satisfies the lower bound \(c\) and upper
bound \(C\) in \eqref{eq:average-embedding-bounds}.
For \(u\in\U(M_j)\) and \(n\ge j\), right invariance of Haar
measure gives \(\|F_n(ux)\|_{L_p(D\bar\otimes N)}
=\|F_nx\|_{L_p(D\bar\otimes N)}\).
Since the free ultrafilter contains \(\{n:n\ge j\}\), we have
\(\|J_0(ux)\|_{L_p(B_0)}=\|J_0x\|_{L_p(B_0)}\) for
\(u\in\bigcup_j\U(M_j)\).
To extend this identity to \(\U(M)\), write \(u=e^{ia}\) with
\(a=a^*\in M\), \(\|a\|_\infty\le\pi\).
Kaplansky's density theorem gives uniformly bounded self-adjoint
elements of \(\bigcup_jM_j\) converging strongly to \(a\).
Their exponentials converge to \(u\) in the strong-\(*\) topology.
Continuity of bounded multiplication on \(L_p(M)\) therefore gives
\begin{equation}\label{eq:first-average-left-invariance}
 \|J_0(ux)\|_{L_p(B_0)}=\|J_0x\|_{L_p(B_0)}
 \qquad(u\in\U(M),\ x\in L_p(M)).
\end{equation}

For the second average, define
\[
 F_n^{\,r}:L_p(M)\longrightarrow L_p(D\bar\otimes B_0),\qquad
 (F_n^{\,r}x)(t)=J_0(xu_n(t)).
\]
The right-sided estimate in \eqref{eq:one-sided-rademacher-bound},
applied to \(J_0\), and the argument of
\eqref{eq:average-rademacher-bound} prove that
\(\{F_n^{\,r}x:n\ge1\}\) is \(p\)-equiintegrable for each \(x\).
The support transfer and separable-subalgebra construction in
Subsection~\ref{subsec:ambient-transfer} yield a properly infinite
von Neumann algebra \(B\) with separable predual and an embedding
\(J:L_p(M)\to L_p(B)\).
Left invariance of Haar measure gives right unitary invariance of
the norm of \(J\). Equation~\eqref{eq:first-average-left-invariance},
applied at each \(t\), preserves left unitary invariance. Thus,
\begin{equation}\label{eq:averaged-biunitary-invariance}
 \begin{gathered}
 c\|x\|_{L_p(M)}\le\|Jx\|_{L_p(B)}\le C\|x\|_{L_p(M)},\\
 \|J(uxv)\|_{L_p(B)}=\|Jx\|_{L_p(B)}
 \qquad(x\in L_p(M),\ u,v\in\U(M)),\\
 \mathcal F(N)\longrightarrow\mathcal F(B).
 \end{gathered}
\end{equation}
Lemma~\ref{lem:twoindex-disjointification} converts \(J\) into an
isometric embedding in a further target algebra whose flow still
receives a map from \(\mathcal F(N)\).

\begin{lemma}[Disjointification of a bi-unitarily invariant embedding]\label{lem:twoindex-disjointification}
Let $1\le p<2$, let $M$ be a properly infinite hyperfinite von Neumann
algebra with separable predual, and let $B$ be a properly infinite
von Neumann algebra with separable predual. Suppose that
$J:L_p(M)\longrightarrow L_p(B)$ is a linear embedding and that
$0<c\le C<\infty$ satisfy
\begin{equation}\label{eq:disj-embedding-bounds}
 c\|x\|_{L_p(M)}\le \rho(x):=\|Jx\|_{L_p(B)}
 \le C\|x\|_{L_p(M)}\qquad(x\in L_p(M)).
\end{equation}
Assume also that
\begin{equation}\label{eq:disj-biunitary-invariance}
 \rho(uxv)=\rho(x)
 \qquad(x\in L_p(M),\ u,v\in\mathcal U(M)).
\end{equation}
There is a linear isometric embedding of $L_p(M)$ into $L_p(B_2)$,
where $B_2$ is obtained from $B$ by two successive Ocneanu ultrapowers,
with passage after each ultrapower to a subalgebra with separable predual
and a faithful normal conditional expectation. Each chosen subalgebra
contains the diagonal copy of the algebra used to form that ultrapower.
\end{lemma}

\begin{proof}
Choose isometries $s_i\in M$, $i\ge1$, with $s_i^*s_j=\delta_{ij}1$.
For elements $\xi,\eta\in L_p(R)$, where $R$ is any von Neumann
algebra, disjointness means $\xi^*\eta=\xi\eta^*=0$.
For $i,j,k\ge1$, define
\begin{equation}\label{eq:disj-array-and-sums}
 \begin{aligned}
 H_{ij}&:L_p(M)\longrightarrow L_p(B),&
 H_{ij}x&:=J(s_i x s_j^*),\\
 a_k&:L_p(M)\longrightarrow[0,\infty),&
 a_k(x)&:=\rho\!\left(\sum_{\ell=1}^k s_\ell x s_\ell^*\right).
 \end{aligned}
\end{equation}
Fix a free ultrafilter $\omega$ on $\mathbb N$. Let $\mathcal V$ be the
product ultrafilter on $\mathbb N^2$ defined by
\[
 \lim_{(i,j)\to\mathcal V}b_{ij}
 :=\lim_{j\to\omega}\lim_{i\to\omega}b_{ij}
 \qquad((b_{ij})\text{ a bounded scalar family}).
\]
Let $G$ be the Groh--Raynaud ultrapower of $B$ for which
$L_p(B)^{\mathcal V}$ identifies isometrically with $L_p(G)$.
Choose a faithful normal state $\psi$ on $B$, and let $e\in G$ be
the support of the constant class of $\psi$ in $L_1(B)^{\mathcal V}$.
Since $B$ is $\sigma$-finite, $e$ is also the support of the diagonal
copy of $L_p(B)$ in $L_p(G)$
\cite[discussion preceding Theorem~4.6]{RX}. Define
\begin{equation}\label{eq:disj-singular-compression}
 \begin{aligned}
 H,D&:L_p(M)\longrightarrow L_p(G),\\
 Hx&:=(H_{ij}x)^{\mathcal V},&
 Dx&:=(1-e)Hx(1-e),\\
 \lambda(x)&:=\|Dx\|_{L_p(G)}&&\qquad(x\in L_p(M)).
 \end{aligned}
\end{equation}

\smallskip
\noindent\emph{Independent copies of the singular part.}
For each $x\in L_p(M)$, Theorem~4.6 of \cite{RX} supplies a bounded
$p$-equiintegrable representative $(g_a^x)_{a\in\mathbb N^2}$ in $L_p(B)$
of the element $Hx-Dx\in L_p(G)$. For $a=(i,j)$, put
$H_ax=H_{ij}x$ and $t_a^x=H_ax-g_a^x$. Then
$(t_a^x)^{\mathcal V}=Dx$, and $Dx$ is disjoint from the diagonal copy
of every $z\in L_p(B)$. Consequently, for each fixed $z\in L_p(B)$,
we have
\begin{equation}\label{eq:disj-cross-products}
 \lim_{a\to\mathcal V}\left(
    \|z^*t_a^x\|_{L_{p/2}(B)}
   +\|z(t_a^x)^*\|_{L_{p/2}(B)}\right)=0.
\end{equation}
The norms in \eqref{eq:disj-cross-products} are the $L_{p/2}(B)$
quasi-norms. Multiplication, adjoints and absolute powers are compatible
with the Groh--Raynaud identifications of the $L_p$-ultrapowers
\cite[Section~1 and Facts~1.1--1.3]{RX}.

For $k\ge1$, let $\mathcal V^{(k)}$ be the product ultrafilter on
$(\mathbb N^2)^k$, with indices $a_\ell=(i_\ell,j_\ell)$, $1\le\ell\le k$.
The limit in $a_k$ is taken first, followed by the limits in
$a_{k-1},\ldots,a_1$. For fixed $x_1,\ldots,x_k\in L_p(M)$,
the $k$ elements represented by $(t_{a_\ell}^{x_\ell})$ are mutually
disjoint in $L_p(B)^{\mathcal V^{(k)}}$.
Indeed, for $\ell<m$ and a fixed value of $a_\ell$, apply
\eqref{eq:disj-cross-products} to $z=t_{a_\ell}^{x_\ell}\in L_p(B)$.
Both cross products vanish in the limit in $a_m$.
Thus, for every $(\varepsilon_\ell)_{\ell=1}^k\in\{-1,1\}^k$, we have
\begin{equation}\label{eq:disj-sum-norm}
 \lim_{\mathcal V^{(k)}}\left\|
            \sum_{\ell=1}^k\varepsilon_\ell t_{a_\ell}^{x_\ell}
                     \right\|_{L_p(B)}
       =\left(\sum_{\ell=1}^k\lambda(x_\ell)^p\right)^{1/p}.
\end{equation}
The row indices $i_1,\ldots,i_k$ are pairwise distinct on a
set belonging to $\mathcal V^{(k)}$, and the same holds for the column
indices $j_1,\ldots,j_k$. This follows because $\omega$ contains the
complement of every finite subset of $\mathbb N$.

\smallskip
\noindent\emph{Identification of the singular norm.}
Fix $x\in L_p(M)$. For every tuple of $k$ pairs
$(i_\ell,j_\ell)$ with distinct row indices and distinct column indices,
we have
\begin{equation}\label{eq:disj-randomized-copies}
 a_k(x)=\mathbb E_\varepsilon\left\|\sum_{\ell=1}^k
                    \varepsilon_\ell H_{i_\ell j_\ell}x\right\|_{L_p(B)}.
\end{equation}
Here $\mathbb E_\varepsilon$ denotes the average over $\{-1,1\}^k$.
To obtain \eqref{eq:disj-randomized-copies}, use
\eqref{eq:disj-biunitary-invariance} with finite permutations of the
range projections $s_i s_i^*$ on the left and on the right, and then
multiply the left range projections by the signs $\varepsilon_\ell$.
Since $(g_a^x)_{a\in\mathbb N^2}$ is $p$-equiintegrable,
\cite[Lemma~5.2]{RX} gives
$\mathbb E_\varepsilon\|\sum_{\ell=1}^k\varepsilon_\ell
g_{a_\ell}^x\|_{L_p(B)}=o(k^{1/p})$ uniformly in all tuples.
Lemma~5.2 of \cite{RX}, stated for sequences, applies to this family
because $\mathbb N^2$ is countable. Substitute
$H_{a_\ell}x=g_{a_\ell}^x+t_{a_\ell}^x$ in
\eqref{eq:disj-randomized-copies} and take the $\mathcal V^{(k)}$-limit.
The finite sign average commutes with the ultralimit. By
\eqref{eq:disj-sum-norm} and the triangle inequality, we obtain
\begin{equation}\label{eq:disj-singular-norm-limit}
 \lim_{k\to\infty} k^{-1/p}a_k(x)=\lambda(x).
\end{equation}
Since the elements $s_\ell x s_\ell^*$, $1\le\ell\le k$, are mutually
disjoint, \eqref{eq:disj-embedding-bounds} and
\eqref{eq:disj-singular-norm-limit} give
\begin{equation}\label{eq:disj-singular-norm-bounds}
 c\|x\|_{L_p(M)}\le\lambda(x)\le C\|x\|_{L_p(M)}
 \qquad(x\in L_p(M)).
\end{equation}
For $u,v\in\mathcal U(M)$, apply $u$ and $v$ in each of the $k$ range
projections in \eqref{eq:disj-array-and-sums}. By
\eqref{eq:disj-biunitary-invariance} and
\eqref{eq:disj-singular-norm-limit}, we have
$\lambda(uxv)=\lambda(x)$ for every $x\in L_p(M)$.

Let $x,y\in L_p(M)_+$ have orthogonal support projections $b,d\in M$.
For each tuple of $2k$ independent indices with all rows distinct and
all columns distinct, we have
\[
 a_k(x+y)=
 \mathbb E_\varepsilon\left\|\sum_{\ell=1}^k
 \varepsilon_\ell H_{i_\ell j_\ell}x+
 \sum_{\ell=1}^k\varepsilon_{k+\ell}H_{i_{k+\ell}j_{k+\ell}}y
 \right\|_{L_p(B)}.
\]
Indeed, choose the left unitary to extend the partial isometry
$\sum_{\ell=1}^k (s_{i_\ell} b s_\ell^*
                 +s_{i_{k+\ell}}d s_\ell^*)$.
The complements of the initial and final projections each contain an
unused range projection $s_m s_m^*$ equivalent to $1$, so both complements
are equivalent to $1$. This permits the extension to a unitary.
Choose the right unitary by the adjoint construction with the $j$ indices,
and insert the signs on the left. The union of the two families
$(g_a^x)_a$ and $(g_a^y)_a$ is $p$-equiintegrable. By
\cite[Lemma~5.2]{RX}, the sign sums of these regular parts have norm
$o(k^{1/p})$ uniformly in the tuples. Applying
\eqref{eq:disj-sum-norm} to the remaining $2k$ singular parts and
\eqref{eq:disj-singular-norm-limit} to $a_k(x+y)$ gives
\begin{equation}\label{eq:disj-orthogonal-additivity}
 \lambda(x+y)^p=\lambda(x)^p+\lambda(y)^p.
\end{equation}
Define
$F:M_*^+\to[0,\infty)$ by $F(\varphi)=\lambda(\varphi^{1/p})^p$.
The function $F$ is homogeneous of degree one and, by
\eqref{eq:disj-orthogonal-additivity}, additive on functionals with
orthogonal supports. By \eqref{eq:disj-singular-norm-bounds}, we have
$F(\varphi)\le C^p\|\varphi\|$ for $\varphi\in M_*^+$.
The norm continuity of $F$ follows from continuity of the positive Mazur
map $\varphi\mapsto\varphi^{1/p}$ from $M_*^+$ to $L_p(M)_+$.
Thus, $F$ is a continuous finite measure in the terminology of
\cite[Definition~4.1]{Sherman}. Since $M$ is hyperfinite
and has no type $\mathrm I_2$ summand, it has $\mathrm{EP}_1$
\cite[Definition~4.1 and Theorem~1.2]{Sherman}.
Consequently, there is $h\in M_+$ such that
$F(\varphi)=\varphi(h)$ for every $\varphi\in M_*^+$.
Equation~\eqref{eq:disj-singular-norm-bounds} gives
$c^p1\le h\le C^p1$. The equality $\lambda(u\xi u^*)=\lambda(\xi)$
implies $h\in Z(M)$. Left unitary invariance of $\lambda$ and the
Russo--Dye theorem give
\[
 \lambda(a\xi)\le\|a\|_\infty\lambda(\xi)
 \qquad(a\in M,\ \xi\in L_p(M)).
\]
For the partial isometry $v\in M$ in the polar decomposition $x=v|x|$,
apply this inequality to $v$ and $v^*$ to obtain
$\lambda(x)=\lambda(|x|)$. The definition of $F$ and the centrality
of $h$ now give
\begin{equation}\label{eq:disj-central-weight-norm}
 \lambda(x)=\|h^{1/p}x\|_{L_p(M)} \qquad(x\in L_p(M)).
\end{equation}

\smallskip
\noindent\emph{Realization in the two Ocneanu ultrapowers.}
The map $D$ in \eqref{eq:disj-singular-compression} takes values in
$L_p(G)$. We now realize its norm $\lambda$ in the $L_p$-space of an
algebra obtained by two successive Ocneanu ultrapowers.
Fix a countable dense sequence $(x_m)_{m\ge1}$ in $L_p(M)$.
Let $d\in G$ be the join of the left and right support projections of
$Dx_m$, $m\ge1$. The projections $d$ and $e$ are $\sigma$-finite,
and $d\le1-e$. By \cite[Theorem~2.3]{RX}, there are projections
$f_{ij}\in B$, $i,j\ge1$, whose coordinatewise action defines a
projection $f\in G$ satisfying $d\le f\le1-e$. Thus,
$fHx_mf=Dx_m$ for every $m$, and by continuity
$fHxf=Dx$ for every $x\in L_p(M)$. Consequently, we have
\begin{equation}\label{eq:disj-iterated-compression-norm}
 \lambda(x)=\lim_{j\to\omega}\lim_{i\to\omega}
                  \|f_{ij}J(s_i x s_j^*)f_{ij}\|_{L_p(B)}
 \qquad(x\in L_p(M)).
\end{equation}
The projections $f_{ij}$ represent the operator $f$ through their
action on $L_p(B)^{\mathcal V}$.

Let $G_1$ be the Groh--Raynaud ultrapower of $B$ in the $i$ variable,
and let $e_1\in G_1$ be the support of the constant class of $\psi$.
Define
\[
 H_1:L_p(M)\longrightarrow L_p(G_1),\qquad
 H_1x:=(J(s_i x))_i^\omega.
\]
For each $j\ge1$, let $b_j\in G_1$ be the projection represented by
the coordinatewise action of $(f_{ij})_{i\ge1}$.
For each fixed $x\in L_p(M)$, the family $(J(s_i x))_{i\ge1}$ is
$p$-equiintegrable by the one-sided estimate
\eqref{eq:one-sided-rademacher-bound} in Subsection~\ref{subsec:averaging}.
Indeed, the factorization $x=x_2x_r$, with $x_2\in L_2(M)$,
$x_r\in L_r(M)$ and $1/p=1/2+1/r$, gives a sign-sum bound of order
$\sqrt{k}$, whereas failure of $p$-equiintegrability would give a
lower bound of order $k^{1/p}$.
By \cite[Theorem~4.6]{RX}, we have
$(1-e_1)H_1x(1-e_1)=0$ for every $x\in L_p(M)$.
It follows that $H_1x$ and $b_jH_1xb_j$, $j\ge1$, are supported by
the central support $z_{G_1}(e_1)$.

Let $r_1\in G_1$ be the join of $e_1$ and the left and right support
projections of $H_1x_m$ and $b_jH_1x_mb_j$, $j,m\ge1$.
Then $r_1$ is $\sigma$-finite and $r_1\le z_{G_1}(e_1)$.
Since $e_1G_1e_1\cong B^\omega$ is properly infinite, the projection
comparison argument in Subsection~\ref{subsec:ambient-transfer} gives
$r_1\precsim e_1$. Choose a partial isometry $v_1\in G_1$ with
$v_1^*v_1=r_1$ and $v_1v_1^*\le e_1$, and put
\begin{equation}\label{eq:disj-first-transfer}
 \begin{aligned}
 K_1&:L_p(M)\longrightarrow L_p(B^\omega),&
 K_1x&:=v_1H_1xv_1^*,\\
 d_j&:=v_1r_1b_jr_1v_1^*\in B^\omega&&\qquad(j\ge1).
 \end{aligned}
\end{equation}
Each $d_j$ is a positive contraction. The choice of $r_1$ and the
coordinatewise module action give
\begin{equation}\label{eq:disj-first-compression-norm}
 \|d_jK_1(z)d_j\|_{L_p(B^\omega)}
       =\lim_{i\to\omega}\|f_{ij}J(s_i z)f_{ij}\|_{L_p(B)}
 \qquad(z\in L_p(M)).
\end{equation}
Choose a subalgebra $B_1\subset B^\omega$ with separable predual
that contains the diagonal $B$, all $d_j$, and coefficients
$a_{m,k}\in B^\omega$ satisfying
$a_{m,k}(\psi^\omega)^{1/p}\to K_1(x_m)$ in $L_p(B^\omega)$ as
$k\to\infty$. Include the rational modular translates of these
generators, as in Subsection~\ref{subsec:ambient-transfer}.
The resulting algebra $B_1$ is invariant under $\sigma^{\psi^\omega}$
and admits a $\psi^\omega$-preserving faithful normal conditional
expectation from $B^\omega$. In particular,
$K_1(L_p(M))\subset L_p(B_1)$. Write
$\psi_1=\psi^\omega|_{B_1}$ for the faithful normal state used in
the second ultrapower.

Form the Groh--Raynaud ultrapower $G_2$ of $B_1$ in the $j$ variable,
and let $e_2\in G_2$ be the support of the constant class of $\psi_1$.
Apply the right multiplication version of
\eqref{eq:one-sided-rademacher-bound} to the bounded linear map
$K_1:L_p(M)\to L_p(B_1)$. For each fixed $x\in L_p(M)$, the family
$(K_1(xs_j^*))_{j\ge1}$ is $p$-equiintegrable. Thus, the element
$(K_1(xs_j^*))_j^\omega\in L_p(G_2)$ is supported by $z_{G_2}(e_2)$.
Multiplication on both sides by the element of $G_2$ represented by
$(d_j)_j$ preserves this central support. Define
\[
 K_2:L_p(M)\longrightarrow L_p(G_2),\qquad
 K_2x:=(d_jK_1(xs_j^*)d_j)_j^\omega.
\]
Equations~\eqref{eq:disj-iterated-compression-norm} and
\eqref{eq:disj-first-compression-norm} give
$\|K_2x\|_{L_p(G_2)}=\lambda(x)$ for every $x\in L_p(M)$.
Let $r_2\in G_2$ be the join of $e_2$ and the left and right support
projections of $K_2x_m$, $m\ge1$. Then $r_2$ is $\sigma$-finite,
$r_2\le z_{G_2}(e_2)$, and $r_2\precsim e_2$.
Choose $v_2\in G_2$ with $v_2^*v_2=r_2$ and $v_2v_2^*\le e_2$.
The map
$x\mapsto v_2K_2xv_2^*$ takes values in
$L_p(e_2G_2e_2)=L_p(B_1^\omega)$.
Use the construction in Subsection~\ref{subsec:ambient-transfer}
to select a subalgebra $B_2\subset B_1^\omega$ with separable predual
and a faithful normal expectation, containing the diagonal $B_1$,
such that $v_2K_2(L_p(M))v_2^*\subset L_p(B_2)$.
Finally, define
\[
 I:L_p(M)\longrightarrow L_p(B_2),\qquad
 I(x):=v_2K_2(h^{-1/p}x)v_2^*.
\]
Since $h$ is boundedly invertible and central,
\eqref{eq:disj-central-weight-norm} gives
$\|I(x)\|_{L_p(B_2)}=\lambda(h^{-1/p}x)=\|x\|_{L_p(M)}$
for every $x\in L_p(M)$.
\end{proof}

\subsection{An isometry produces a map of flows}\label{subsec:disj-isometry-to-flow}
Apply Lemma~\ref{lem:twoindex-disjointification} to the map \(J\)
constructed in Subsection~\ref{subsec:averaging}. We obtain a linear
isometry \(I:L_p(M)\to L_p(W)\), where \(W\) has separable predual.
The operations in Subsection~\ref{subsec:ambient-transfer} also give
an equivariant normal unital positive map
\(K_{N,W}:\Zc(N)\to\Zc(W)\).
Since \(M\) is hyperfinite and has no \(\mathrm I_2\) summand,
\cite[Theorem~1.2]{Sherman} supplies a normal Jordan
\(*\)-monomorphism \(j:M\to W\) and a normal positive projection
\(P:W\to j(M)\). Define
\[
 Q:W\longrightarrow M,\qquad Q:=j^{-1}\circ P.
\]
Sherman's formula for \(I(\varphi^{1/p})\) and the isometry of \(I\) give
\[
 (\varphi\circ Q)(1_W)
   =\|I(\varphi^{1/p})\|_{L_p(W)}^p
   =\|\varphi^{1/p}\|_{L_p(M)}^p
   =\varphi(1_M)
 \qquad(\varphi\in M_*^+).
\]
Thus, \(Q(1_W)=1_M\), so \(Q\) is contractive. Since
\(Q\circ j=\mathrm{id}_M\) and \(j\) is contractive, the preadjoint
\[
 Q_*:M_*\longrightarrow W_*,\qquad Q_*(\varphi)=\varphi\circ Q,
\]
is an isometric embedding.

Apply the Kirchberg lemma in the formulation of
Oikhberg--Rosenthal--St{\o}rmer \cite[Lemma~2]{ORS} to \(Q_*\).
There are projections \(r_1,r_2\in W\), a normal \(*\)-homomorphism
\(\alpha:M\to r_1Wr_1\), and a normal \(*\)-antihomomorphism
\(\beta:M\to r_2Wr_2\), such that
\(\ker\alpha\cap\ker\beta=\{0\}\). Each of the two image algebras
admits a faithful normal conditional expectation from its respective
corner. Let \(z\in Z(M)\) satisfy \(\ker\alpha=(1-z)M\).
Then \(\alpha|_{zM}\) and \(\beta|_{(1-z)M}\) are faithful.
Compress the expectation onto \(\beta(M)\) by the projection
\(\beta(1-z)\). For each nonzero summand, the corner and expectation
constructions in Subsection~\ref{subsec:ambient-transfer} give
equivariant normal unital positive maps
\[
 \Zc(W)\longrightarrow\Zc(zM),\qquad
 \Zc(W)\longrightarrow\Zc(((1-z)M)^{\rm op}).
\]

To identify the last codomain with \(\Zc((1-z)M)\), let \(S\) be
a \(\sigma\)-finite von Neumann algebra and let \(\gamma\) be a
faithful normal state on \(S\). In the core isomorphism
\(\core(S^{\rm op})\cong\core(S)^{\rm op}\), the element
\(a^{\rm op}\), \(a\in S\), is identified with the corresponding
element of \(S^{\rm op}\subset\core(S)^{\rm op}\), and
\(\lambda_{\gamma^{\rm op}}(t)\) is sent to
\(\lambda_\gamma(t)^{\rm op}\). The covariance relations agree because
\[
 \sigma_t^{\gamma^{\rm op}}(a^{\rm op})
   =(\sigma_{-t}^\gamma(a))^{\rm op}
 \qquad(a\in S,\ t\in\mathbb R).
\]
Both dual actions multiply the corresponding implementing unitaries
by \(e^{-ist}\), for \(s,t\in\mathbb R\). Therefore, the induced
identification \(\Zc(S^{\rm op})\cong\Zc(S)\) intertwines the dual
actions with the same time parameter.
Apply this identification to \(S=(1-z)M\). Combining the two
component maps gives an equivariant normal unital positive map
\(K_{W,M}:\Zc(W)\to\Zc(M)\). Its composition with \(K_{N,W}\) gives
\[
             \mathcal F(N)\longrightarrow\mathcal F(W)
                         \longrightarrow\mathcal F(M).
\]
Thus, the Banach embedding \(T:L_p(M)\to L_p(N)\) yields
\(\mathcal F(N)\to\mathcal F(M)\), as required for the necessary
condition in Theorem~\ref{thm:criterion}(1).

\section{Realizing a map of flows inside the original target}

\subsection{A modularly invariant finite-dimensional approximation}

\begin{lemma}\label{lem:invariant-filtration}
Let $P$ be a hyperfinite type III von Neumann algebra with separable
predual. There exist a faithful normal state $\chi'$ on $P$ and an
increasing sequence of unital finite-dimensional subalgebras $D_n\subset P$
such that
\[
 \sigma_t^{\chi'}(D_n)=D_n\qquad(n\ge1,\ t\in\mathbb R),
 \qquad
 \overline{\bigcup_{n\ge1}D_n}^{\,\mathrm{weak}^*}=P.
\]
\end{lemma}

\begin{proof}
The proof of \cite[Theorem~7.1, pp.~87--92]{JRX} gives a central
decomposition
\[
 P=\int_X^\oplus P_x\,d\mu(x),\qquad
 Z(P)=L_\infty(X,\mu),
\]
where $(X,\mu)$ is a standard probability space. For almost every
$x\in X$, that proof also gives a faithful normal state $\chi_x$ on
$P_x$ and increasing unital finite-dimensional subalgebras
$A_n(x)\subset P_x$ with $\chi_x$-preserving faithful normal conditional
expectations and weak$^*$-dense union. The states and subalgebras form
measurable fields. Put
\[
 \chi(b)=\int_X\chi_x(b_x)\,d\mu(x)
 \quad\left(b=\int_X^\oplus b_x\,d\mu(x)\in P\right),\qquad
 \mathcal A_n=\int_X^\oplus A_n(x)\,d\mu(x),\qquad
 A_0(x)=\mathbb C1.
\]
The decomposable algebras $\mathcal A_n$ are type I and generate $P$.
We modify $\chi$ and then choose finite-dimensional subalgebras of
the $\mathcal A_n$.

For the inclusions $A_n(x)\subset A_{n+1}(x)$, choose compatible
matrix units indexed by paths in the Bratteli diagram of the increasing
sequence $(A_n(x))_n$. Write $e^n_{\alpha\beta}(x)$ for the matrix
unit indexed by two length-$n$ paths $\alpha,\beta$ ending at the
same block $i$ of $A_n(x)$. Compatibility means that
\begin{equation}\label{eq:filtration-matrix-inclusion}
 e^n_{\alpha\beta}(x)
 =\sum_{\substack{e\text{ an edge}\\
                   \text{from block }i}}
     e^{n+1}_{\alpha e,\beta e}(x).
\end{equation}
Parallel edges are counted separately. The choices of matrix units
can be made measurably. To see this, partition $X$ at each level by the
countably many possible block dimensions and inclusion multiplicities.
On each part, the matrix-unit equations, the identity decomposition
and \eqref{eq:filtration-matrix-inclusion} are Borel equations for
finitely many unit-ball sections, with respect to the strong$^*$
Borel structure. Require in addition that $\chi_x$ vanish on the
off-diagonal matrix units. The relation defined by these equations
is nonempty, as the density computation in the next paragraph shows.
A universally measurable selection is measurable for the completed
measure $\mu$. Discarding the union of the exceptional null sets over
all levels gives compatible choices for almost every $x$.

Use the unnormalized matrix trace on each block to represent the
restriction of $\chi_x$. Denote its density in block $i$ of
$A_n(x)$ by $D_{n,i}(x)$. Modular invariance of $A_n(x)$ implies that
the density in block $j$ of $A_{n+1}(x)$ has the form
\[
 D_{n+1,j}(x)
 =\bigoplus_i D_{n,i}(x)\otimes b^n_{ij}(x),
\]
where $b^n_{ij}(x)$ is a positive invertible matrix on the
multiplicity space for the inclusion of block $i$ into block $j$.
Indeed, on the $i$-component of block $j$, differentiation of the
equality of the two modular groups shows that
$\log D_{n+1,j}(x)-\log D_{n,i}(x)\otimes1$ commutes with the
embedded block $i$. Diagonalize each $b^n_{ij}(x)$ on its
multiplicity space. Its eigenvalues give positive weights $w_e(x)$
for the edges $e$ from level $n$ to level $n+1$. Consistency of the
restricted states gives
\begin{equation}\label{eq:filtration-edge-weights}
 \sum_{\substack{e\text{ an edge}\\
                   \text{from block }i}}w_e(x)=1,\qquad
 a_\alpha(x):=\chi_x(e^n_{\alpha\alpha}(x))
             =\prod_{e\in\alpha}w_e(x)
 \quad(|\alpha|=n).
\end{equation}
In particular, the densities are diagonal in the path matrix units
at every level. The weights $w_e(x)$ are measurable, since they are
ratios of the measurable diagonal values $a_\alpha(x)$.

For each block at level $n-1$, choose positive rational weights
$(w'_e(x))_e$ on its outgoing edges such that
\begin{equation}\label{eq:filtration-rational-weights}
 \sum_e w'_e(x)=1,\qquad
 \left|\log\frac{w'_e(x)}{w_e(x)}\right|<2^{-n}.
\end{equation}
For each possible number of outgoing edges, enumerate the positive
rational probability vectors and choose the first vector satisfying
\eqref{eq:filtration-rational-weights}. The choice is measurable.
Put $a'_\alpha(x)=\prod_{e\in\alpha}w'_e(x)$ and define $h_n\in\mathcal A_n$
by
\begin{equation}\label{eq:filtration-perturbation}
 (h_n)_x=\sum_{|\alpha|=n}
          \frac{a'_\alpha(x)}{a_\alpha(x)}e^n_{\alpha\alpha}(x).
\end{equation}
The operators $h_n$ commute, belong to $P_\chi$, and satisfy
\[
 \|\log h_{n+1}-\log h_n\|_P\le2^{-(n+1)}.
\]
Consequently, $(h_n)_n$ converges in operator norm to a positive
boundedly invertible $h\in P_\chi$. For $m\ge n$, the operator
$h_mh_n^{-1}$ commutes with $\mathcal A_n$. In fact, its diagonal
coefficient on a length-$m$ path is the product of the ratios
$w'_e(x)/w_e(x)$ on the edges after level $n$. That coefficient is
unchanged when one length-$n$ prefix is replaced by another ending
at the same block. Passing to the norm limit shows that $hh_n^{-1}$
also commutes with $\mathcal A_n$.

Define a faithful normal state on $P$ by
\[
 \chi'(b)=\frac{\chi(h^{1/2}bh^{1/2})}{\chi(h)}
 \qquad(b\in P).
\]
The modular group of this centralizer perturbation satisfies
$\sigma_t^{\chi'}(b)=h^{it}\sigma_t^\chi(b)h^{-it}$ for
$b\in P$ and $t\in\mathbb R$. Combining this formula with
\eqref{eq:filtration-perturbation} and the commutation of $hh_n^{-1}$
with $\mathcal A_n$, we obtain
\begin{equation}\label{eq:filtration-modular-action}
 \sigma_t^{\chi'}(e^n_{\alpha\beta})(x)
 =\left(\frac{a'_\alpha(x)}{a'_\beta(x)}\right)^{it}
     e^n_{\alpha\beta}(x)
 \qquad(t\in\mathbb R,\ |\alpha|=|\beta|=n).
\end{equation}
The formula is understood on each measurable part of $X$ where the
path labels $\alpha,\beta$ are defined.

For each $n$, let $\mathcal P_n$ be the countable measurable partition
of $X$ determined by the labelled diagrams and the rational weights
$w'_e$ through level $n$, together with the first $n$ sets in a
fixed countable generator of the measure algebra of $(X,\mu)$.
Let $B_n\subset\mathcal A_n$ consist of the sections whose matrix
coefficients in the chosen path matrix units are constant on each
cell of $\mathcal P_n$. Thus, $B_n$ is a countable product of
finite-dimensional matrix algebras. The partitions refine as $n$
increases, and \eqref{eq:filtration-matrix-inclusion} gives
$B_n\subset B_{n+1}$. Equation~\eqref{eq:filtration-modular-action}
gives $\sigma_t^{\chi'}(B_n)=B_n$ for every $t\in\mathbb R$.
Moreover, $\bigcup_n B_n$ generates $P$. Indeed, its weak$^*$-closure
contains $Z(P)$ by the choice of the generating sets of the measure
algebra. Multiplying the matrix units that are constant on each cell by elements
of $Z(P)$ then gives every decomposable algebra $\mathcal A_n$.

Choose a union $G_n$ of finitely many cells of $\mathcal P_n$ with
$\mu(X\setminus G_n)<2^{-n}$. Then $I_n=1_{G_n}B_n$ is a
finite-dimensional two-sided ideal in $B_n$. Set
\[
 D_0=\mathbb C1,\qquad D_n=D_{n-1}+I_n\quad(n\ge1).
\]
Since $D_{n-1}\subset B_n$ and $I_n$ is an ideal in $B_n$,
the algebra $D_n$ is a unital finite-dimensional $*$-subalgebra of
$B_n$. Both $D_{n-1}$ and $I_n$ are invariant under
$\sigma^{\chi'}$, so $D_n$ is invariant as well.
For $b\in B_k$ and $n\ge k$, we have $1_{G_n}b\in I_n\subset D_n$.
The summability of $\mu(X\setminus G_n)$ implies that
$1_{G_n}\to1$ strongly, and hence $1_{G_n}b\to b$ strongly.
Thus, the weak$^*$-closure of $\bigcup_nD_n$ contains every $B_k$
and therefore equals $P$.
\end{proof}

\subsection{From an equivariant Markov map to an inclusion with expectation}

For \(R=M,N\), let \(\tau_R\) be the canonical trace on
\(\core(R)\). We use the convention
\(\tau_R\circ\theta_t^R=e^{-t}\tau_R\), \(t\in\mathbb R\), for the
dual action. Equivariant Markov maps are defined in
\eqref{eq:flow-map-definition}.

\begin{proposition}\label{prop:markov-to-expected}
Let \(M\) and \(N\) be nonzero \(\sigma\)-finite von Neumann algebras
of pure type III. Suppose that \(M\) admits a faithful normal state
\(\varphi\) and an increasing sequence of unital finite-dimensional
\(\sigma^\varphi\)-invariant subalgebras \(A_n\subset M\) whose union
is weakly dense in \(M\).
If there exists an equivariant Markov map
\[
 K:\Zc(N)\longrightarrow\Zc(M),
\]
then there exist a nonzero projection \(q\in N\), a normal unital
\(*\)-monomorphism \(j:M\longrightarrow qNq\), and a faithful normal
conditional expectation \(E:qNq\longrightarrow j(M)\).
Consequently, \(L_p(M)\) is isometric to a \(1\)-complemented subspace
of \(L_p(N)\) for every \(1\le p<\infty\).
If \(N\) is a factor, then \(qNq\) is isomorphic to \(N\).
\end{proposition}

Lemma~\ref{lem:invariant-filtration} supplies the state \(\varphi\) and
the subalgebras \(A_n\) when \(M\) is hyperfinite with separable predual.
The algebras \(A_n\) may have more than one matrix summand.

\begin{proof}
For \(R=M,N\) and \(\beta\in R_*^+\), let \(h_\beta\) be the
Haagerup density of \(\beta\), viewed as an operator affiliated with
\(\core(R)\). Define the positive normal functional \(\widehat\beta\)
on \(\Zc(R)\) by
\begin{equation}\label{eq:expected-spectral-functional}
 \widehat\beta(z)=\tau_R\bigl(z\,1_{(1,\infty)}(h_\beta)\bigr)
 \qquad(z\in\Zc(R)).
\end{equation}
Haagerup--St{\o}rmer \cite[Main Theorem and Proposition~10.1]{HS}
show that
\(\|\widehat\beta\|_{\Zc(R)_*}=\beta(1)\), and that
\(\widehat\beta=\widehat\gamma\) if and only if the norm closures in
\(R_*\) of the unitary orbits of \(\beta,\gamma\in R_*^+\) agree.
Since \(R\) is properly infinite and has no type I summand, the range of
\(\beta\mapsto\widehat\beta\) is
\begin{equation}\label{eq:expected-spectral-range}
 \left\{\nu\in\Zc(R)_*^+:
       \nu\circ\theta_t^R\ge e^{-t}\nu\quad(t\ge0)\right\}.
\end{equation}
The direction of the inequality in
\eqref{eq:expected-spectral-range} follows also from
\[
 \widehat\beta\circ\theta_t^R(z)
 =e^{-t}\tau_R\bigl(z\,1_{(e^{-t},\infty)}(h_\beta)\bigr)
 \ge e^{-t}\widehat\beta(z)
 \qquad(z\in\Zc(R)_+,\ t\ge0).
\]

For \(w>0\), we have
\(1_{(1,\infty)}(w h_\beta)
=\theta_{-\log w}^R(1_{(1,\infty)}(h_\beta))\).
Using \eqref{eq:expected-spectral-functional} and the scaling of
\(\tau_R\), we obtain
\begin{equation}\label{eq:expected-spectral-scaling}
 \widehat{w\beta}
   =w\,\widehat\beta\circ\theta_{\log w}^R
 \qquad(\beta\in R_*^+,\ w>0).
\end{equation}
Spectral calculus on orthogonal corners shows that
\(\widehat{\beta+\gamma}=\widehat\beta+\widehat\gamma\) whenever
\(\beta,\gamma\in R_*^+\) have orthogonal supports. For an isometry
\(v\in R\), define \(v\beta v^*\in R_*^+\) by
\((v\beta v^*)(x)=\beta(v^*xv)\), \(x\in R\).
Its density is \(v h_\beta v^*\). Traciality of \(\tau_R\) therefore gives
\begin{equation}\label{eq:expected-spectral-isometry}
 \widehat{v\beta v^*}=\widehat\beta.
\end{equation}

The preadjoint
\[
 K_*:\Zc(M)_*\longrightarrow\Zc(N)_*,\qquad K_*(\nu)=\nu\circ K,
\]
is positive, preserves the norms of positive normal functionals, and
satisfies
\begin{equation}\label{eq:expected-preadjoint-equivariance}
 K_*(\nu\circ\theta_t^M)=(K_*\nu)\circ\theta_t^N
 \qquad(\nu\in\Zc(M)_*,\ t\in\mathbb R).
\end{equation}
It follows from \eqref{eq:expected-spectral-range} that
\(K_*(\widehat\beta)\) is the spectral functional of a normal state on
\(N\) whenever \(\beta\) is a normal state on \(M\).

We next describe the restricted states at two successive levels
\(A_n\subset A_{n+1}\). In the formulas for these two levels, we suppress
the index \(n\) and use a prime for quantities at level \(n+1\).
Write
\[
 A_n=\bigoplus_{i\in I_n}\mathbb M_{r_i},
 \qquad
 A_{n+1}=\bigoplus_{j\in I_{n+1}}\mathbb M_{r'_j},
\]
where \(I_n,I_{n+1}\) are finite index sets. Choose matrix units
\((e^i_{ab})_{a,b=1}^{r_i}\) so that the density of
\(\varphi|_{A_n}\), with respect to the unnormalized matrix trace
on each summand, is
\[
 D_i=\operatorname{diag}(d_{i,a})_{a=1}^{r_i},
 \qquad d_{i,a}>0,\qquad \sum_{i\in I_n}\sum_{a=1}^{r_i}d_{i,a}=1.
\]
Let \(m_{ij}\) be the multiplicity of the \(i\)-th summand of \(A_n\)
in the \(j\)-th summand of \(A_{n+1}\).
In the latter summand, index the matrix coordinates by
\((i,a,l)\), where \(1\le a\le r_i\) and \(1\le l\le m_{ij}\).
All sums involving \((i,j,l)\) below include only indices with
\(m_{ij}>0\). We may choose the bases of the multiplicity spaces so that
the density of \(\varphi|_{A_{n+1}}\) on summand \(j\) is
\begin{equation}\label{eq:expected-density-recursion}
 \begin{gathered}
 D'_j=\bigoplus_{i:\,m_{ij}>0}D_i\otimes
          \operatorname{diag}(w_{ij1},\ldots,w_{ijm_{ij}}),
 \qquad w_{ijl}>0,\\
 \sum_{\substack{j\in I_{n+1}\\1\le l\le m_{ij}}}w_{ijl}=1
 \qquad(i\in I_n).
 \end{gathered}
\end{equation}
Indeed, the modular groups of \(\varphi|_{A_n}\) and
\(\varphi|_{A_{n+1}}\) are the restrictions of \(\sigma^\varphi\).
In block \(j\), the action \(\operatorname{Ad}((D'_j)^{it})\)
fixes the image of each central projection of \(A_n\).
On the subspace corresponding to summand \(i\), this action agrees
with \(\operatorname{Ad}(D_i^{it}\otimes1)\) on
\(\mathbb M_{r_i}\otimes1\). Differentiating at \(t=0\) shows that
the difference between the restriction of \(\log D'_j\) and
\((\log D_i)\otimes1\) belongs to \(1\otimes\mathbb M_{m_{ij}}\).
Diagonalizing that difference on the multiplicity space gives the
first identity in \eqref{eq:expected-density-recursion}.
Restricting the resulting state back to \(A_n\) gives the sum of the
weights in the second identity. The matrix units can therefore be
chosen recursively through all levels.

For \(i\in I_n\) and \(1\le a\le r_i\), define normal states on \(M\) by
\begin{equation}\label{eq:expected-corner-states}
 \beta_{i,a}(x)=d_{i,a}^{-1}\varphi(e^i_{aa}x e^i_{aa}),
 \qquad \beta_i=\beta_{i,1}\qquad(x\in M).
\end{equation}
The matrix units satisfy
\(\sigma_t^\varphi(e^i_{ab})=(d_{i,a}/d_{i,b})^{it}e^i_{ab}\).
The modular identity for \(\varphi\) consequently gives
\[
 \beta_{i,a}(x)=\beta_{i,b}(e^i_{ba}x e^i_{ab})
 \qquad(x\in M,\ 1\le a,b\le r_i).
\]
Extend \(e^i_{ab}\) to the unitary in \(A_n\) that interchanges the
\(a\)-th and \(b\)-th coordinates of summand \(i\).
Thus, the states \(\beta_{i,a}\) are unitarily conjugate as \(a\) varies,
and their spectral functionals all equal \(\widehat\beta_i\).
Define \(\beta'_{j,k}\) and \(\beta'_j\) at level \(n+1\) by
\eqref{eq:expected-corner-states} with the primed matrix units and
density coefficients. The same modular calculation gives
\(\widehat{\beta'_{j,k}}=\widehat{\beta'_j}\).

Each diagonal projection in \(A_{n+1}\) belongs to \(M_\varphi\).
For distinct diagonal projections \(e',f'\in A_{n+1}\), we therefore have
\(\varphi(e'xf')=\varphi(f'e'x)=0\), \(x\in M\).
Using the expansion
\(e^i_{11}=\sum_{j,l}e'^j_{(i,1,l),(i,1,l)}\) and
\eqref{eq:expected-density-recursion}, we obtain
\[
 \beta_i=\sum_{j,l}w_{ijl}\,\beta'_{j,(i,1,l)}
 \quad\text{as normal functionals on }M.
\]
The states in the sum have pairwise orthogonal supports.
Orthogonal additivity and \eqref{eq:expected-spectral-scaling} now give
\begin{equation}\label{eq:expected-spectral-recursion}
 \widehat\beta_i
 =\sum_{j,l}w_{ijl}\,
            \widehat{\beta'_j}\circ\theta^M_{\log w_{ijl}}
 \qquad(i\in I_n).
\end{equation}

At each level \(n\), choose a normal state \(\psi_i\) on \(N\) for
each \(i\in I_n\) such that
\(\widehat{\psi_i}=K_*(\widehat\beta_i)\).
Such states exist by \eqref{eq:expected-spectral-range} and
\eqref{eq:expected-preadjoint-equivariance}.
Continue to use primes for the corresponding states at level \(n+1\).
Fix positive numbers \(\varepsilon_n\) with
\(\sum_{n=1}^\infty\varepsilon_n<\infty\).
For each \(i\in I_n\), proper infiniteness of \(N\) supplies isometries
\(\widetilde v_{ijl}\in N\) with pairwise orthogonal ranges summing to \(1\),
indexed by \((j,l)\). The normal state
\(\sum_{j,l}w_{ijl}\widetilde v_{ijl}\psi'_j\widetilde v_{ijl}^*\)
on \(N\) has spectral functional
\[
 \begin{aligned}
 \sum_{j,l}w_{ijl}\widehat{\psi'_j}\circ\theta^N_{\log w_{ijl}}
 &=K_*\left(\sum_{j,l}w_{ijl}
                \widehat{\beta'_j}\circ\theta^M_{\log w_{ijl}}\right)\\
 &\overset{\eqref{eq:expected-spectral-recursion}}{=}
 K_*(\widehat\beta_i)=\widehat{\psi_i},
 \end{aligned}
\]
by \eqref{eq:expected-spectral-scaling}--\eqref{eq:expected-preadjoint-equivariance}.
The Haagerup--St{\o}rmer unitary-orbit criterion therefore gives a unitary
\(u_i\in N\) such that \(v_{ijl}:=u_i\widetilde v_{ijl}\) satisfy
\begin{equation}\label{eq:expected-state-approximation}
 \left\|\psi_i-\sum_{j,l}w_{ijl}v_{ijl}\psi'_j v_{ijl}^*
       \right\|_{N_*}<\varepsilon_n
 \qquad(i\in I_n).
\end{equation}
The isometries \(v_{ijl}\) still have pairwise orthogonal ranges summing
to \(1\) for each fixed \(i\).

We construct compatible representations of \(A_n\) inside \(N\).
At the first level, choose isometries \(s_{i,a}\in N\), indexed by
\(i\in I_1\) and \(1\le a\le r_i\), whose ranges are pairwise orthogonal
and sum to \(1\). Represent \(e^i_{ab}\) by \(s_{i,a}s_{i,b}^*\).
Given the isometries \(s_{i,a}\) at level \(n\), define the isometries
at level \(n+1\) by
\begin{equation}\label{eq:expected-isometry-recursion}
 s'_{j,(i,a,l)}=s_{i,a}v_{ijl}.
\end{equation}
For fixed \(j\), these isometries represent the matrix units of the
\(j\)-th block by
\(e'^j_{k k'}\mapsto s'_{j,k}(s'_{j,k'})^*\).
The ranges at level \(n+1\) are pairwise orthogonal and sum to \(1\).
Moreover,
\[
 \sum_{j,l}s'_{j,(i,a,l)}(s'_{j,(i,b,l)})^*
 =s_{i,a}\left(\sum_{j,l}v_{ijl}v_{ijl}^*\right)s_{i,b}^*
 =s_{i,a}s_{i,b}^*.
\]
The level-\((n+1)\) representation therefore extends the level-\(n\)
representation. Let
\(\pi:\bigcup_n A_n\longrightarrow N\) denote the resulting unital
\(*\)-homomorphism.

At level \(n\), define \(\Phi_n\in N_*^+\) by
\begin{equation}\label{eq:expected-approximating-state}
 \Phi_n(x)=\sum_{i\in I_n}\sum_{a=1}^{r_i}
            d_{i,a}\psi_i(s_{i,a}^*x s_{i,a})
 \qquad(x\in N).
\end{equation}
Then \(\Phi_n\) is a normal state and
\(\Phi_n(\pi(a))=\varphi(a)\) for every \(a\in A_n\).
Combining \eqref{eq:expected-density-recursion},
\eqref{eq:expected-isometry-recursion} and
\eqref{eq:expected-state-approximation}, we obtain
\[
 \|\Phi_{n+1}-\Phi_n\|_{N_*}
 \le \sum_{i,a}d_{i,a}
       \left\|\sum_{j,l}w_{ijl}v_{ijl}\psi'_j v_{ijl}^*-\psi_i
       \right\|_{N_*}
 \le\varepsilon_n.
\]
Thus, \((\Phi_n)\) converges in \(N_*\) to a normal state
\(\Phi\). Set \(q=s(\Phi)\in\operatorname{Proj}(N)\).

Fix a level \(n\). The symbols \(s_{i,a},d_{i,a}\) in the rest of this
paragraph refer to this level. Iterating
\eqref{eq:expected-density-recursion} and
\eqref{eq:expected-isometry-recursion} shows that, for every \(m\ge n\),
there are normal states \(\eta_{m,i}\) on \(N\) such that
\[
 \Phi_m(x)=\sum_{i\in I_n}\sum_{a=1}^{r_i}
       d_{i,a}\eta_{m,i}(s_{i,a}^*x s_{i,a})
 \qquad(x\in N).
\]
For fixed \(i\), the state \(\eta_{m,i}\) is independent of \(a\).
For \(x\in N\), we have
\(\|\sum_a s_{i,a}x s_{i,a}^*\|_N=\|x\|_N\), and
\[
 (\Phi_m-\Phi_k)\left(\sum_a s_{i,a}x s_{i,a}^*\right)
 =\left(\sum_a d_{i,a}\right)(\eta_{m,i}-\eta_{k,i})(x)
 \qquad(m,k\ge n).
\]
Consequently,
\[
 \left(\sum_a d_{i,a}\right)
       \|\eta_{m,i}-\eta_{k,i}\|_{N_*}
 \le\|\Phi_m-\Phi_k\|_{N_*}.
\]
Since \(\sum_a d_{i,a}>0\), each sequence \((\eta_{m,i})_{m\ge n}\)
converges in \(N_*\) to a normal state \(\eta_i\).
Passing to the limit in the formula for \(\Phi_m\), we obtain
\begin{equation}\label{eq:expected-limit-state-blocks}
 \Phi(x)=\sum_{i\in I_n}\sum_{a=1}^{r_i}
       d_{i,a}\eta_i(s_{i,a}^*x s_{i,a})
 \qquad(x\in N).
\end{equation}
In particular,
\[
 q=\sum_{i\in I_n}\sum_{a=1}^{r_i}s_{i,a}s(\eta_i)s_{i,a}^*,
\]
so \(q\) commutes with \(\pi(A_n)\).
The restriction of \(\Phi\) to \(qNq\) is faithful.
On the \(i\)-th block of \(q\pi(A_n)\), formula
\eqref{eq:expected-limit-state-blocks} expresses \(\Phi\) as the product
of the positive matrix functional with density \(D_i\) and the faithful restriction of
\(\eta_i\) to its support. Hence, the modular group of \(\Phi|_{qNq}\)
preserves \(q\pi(A_n)\) for every \(n\). It also preserves the von Neumann
algebra
\[
 P=\left(\bigcup_{n\ge1}q\pi(A_n)\right)''\subset qNq.
\]
Takesaki's conditional-expectation theorem gives a faithful normal
conditional expectation \(E:qNq\longrightarrow P\) preserving
\(\Phi|_{qNq}\).

Finally, \eqref{eq:expected-approximating-state} implies
\(\Phi(q\pi(a))=\varphi(a)\) for every \(a\in\bigcup_n A_n\).
The map \(a\mapsto q\pi(a)\) therefore induces a unitary between the
GNS Hilbert spaces of \((M,\varphi)\) and \((P,\Phi|_P)\).
Both states are faithful, so this unitary implements a normal
state-preserving \(*\)-isomorphism \(j:M\longrightarrow P\).
The conditional expectation \(E:qNq\to j(M)\) induces an isometric
embedding of \(L_p(M)\) into \(L_p(qNq)\) with a contractive left inverse.
The corner inclusion \(L_p(qNq)\subset L_p(N)\) is isometric and
has the contractive left inverse \(x\mapsto qxq\).
Composing these two embeddings and their left inverses proves the
\(1\)-complemented assertion for every \(1\le p<\infty\).
\end{proof}

\section{Banach isomorphism and the number of classes}
\begin{proof}[Completion of the proof of Theorem~\ref{thm:criterion}]
Lemma~\ref{lem:invariant-filtration} supplies a faithful normal state
on \(M\) and a modularly invariant finite-dimensional filtration.
Thus, Proposition~\ref{prop:markov-to-expected} applies whenever
\(\mathcal F(N)\to\mathcal F(M)\), and gives a \(1\)-complemented
isometric embedding \(L_p(M)\to L_p(N)\) for every \(1\le p<\infty\).
Together with the necessity proved in Section~2, this proves part~(1).

Suppose that \(\mathcal F(N)\to\mathcal F(M)\) and
\(\mathcal F(M)\to\mathcal F(N)\). Fix \(1\le p<\infty\), and put
\(X=L_p(M)\), \(Y=L_p(N)\). Proposition~\ref{prop:markov-to-expected}
gives complemented embeddings of \(X\) into \(Y\) and of \(Y\) into
\(X\). Moreover, \(X\simeq X\oplus_pX\). Indeed, proper infiniteness
of \(M\) gives isometries \(v_1,v_2\in M\) with
\[
 v_i^*v_j=\delta_{ij}1_M,\qquad v_1v_1^*+v_2v_2^*=1_M.
\]
The bounded linear maps
\[
 \begin{aligned}
 S:X\oplus_pX&\longrightarrow X,&
 S(x_1,x_2)&=v_1x_1+v_2x_2,\\
 R:X&\longrightarrow X\oplus_pX,&
 R(x)&=(v_1^*x,v_2^*x)
 \end{aligned}
\]
satisfy \(RS=I_{X\oplus_pX}\) and \(SR=I_X\).
The same construction in \(N\) proves \(Y\simeq Y\oplus_pY\).
Since the two embeddings are complemented, there are Banach spaces
\(V,W\) with \(Y\simeq X\oplus_pV\) and \(X\simeq Y\oplus_pW\).
Consequently,
\[
 \begin{aligned}
 X\oplus_pY&\simeq X\oplus_pX\oplus_pV
                   \simeq X\oplus_pV\simeq Y,\\
 X\oplus_pY&\simeq Y\oplus_pY\oplus_pW
                   \simeq Y\oplus_pW\simeq X.
 \end{aligned}
\]
Thus, \(X\simeq Y\).

Conversely, if \(L_p(M)\simeq L_p(N)\) for \(1\le p<2\), apply
part~(1) to the isomorphism and its inverse to obtain both maps of
flows. If \(2<p<\infty\), the adjoint of an isomorphism
\(T:L_p(M)\to L_p(N)\) is an isomorphism
\[
 T^*:L_{p'}(N)\longrightarrow L_{p'}(M),\qquad
 p'=\frac{p}{p-1}\in(1,2).
\]
Part~(1), applied to \(T^*\) and \((T^*)^{-1}\), again gives both
maps of flows. This proves part~(2).
\end{proof}

\begin{corollary}\label{cor:continuum-classes}
For each fixed \(1\le p<\infty,\ p\ne2\), the cardinality of the set of
Banach isomorphism classes of the \(L_p\)-spaces in
Theorem~\ref{thm:criterion} is \(2^{\aleph_0}\), even when only factors
of type \(\mathrm{III}_0\) are considered.
\end{corollary}
\begin{proof}
For \(0\le t<2\), let \(M^{(t)}\) be the hyperfinite
type \(\mathrm{III}_0\) factor constructed by Haagerup and Musat
\cite[Theorems~4.1--4.2]{HM}, and write
\(\mathcal F(M^{(t)})=(A_t,\theta^{(t)})\). Their construction gives
\begin{equation}\label{eq:hm-state-displacement}
 \lim_{n\to\infty}
 \|\mu\circ\theta^{(t)}_{2^n}-\mu\|_{(A_t)_*}=t
 \qquad(\mu\text{ a normal state on }A_t).
\end{equation}
Suppose that \(K:A_t\to A_s\) is a normal unital positive linear
map with \(K\circ\theta^{(t)}_r=\theta^{(s)}_r\circ K\) for every
\(r\in\mathbb R\). For any normal state \(\mu\) on \(A_s\),
the functional \(\mu\circ K\) is a normal state on \(A_t\).
Since \(\|K\|=1\), we have
\begin{equation}\label{eq:flow-state-contraction}
\begin{aligned}
 \|(\mu\circ K)\circ\theta^{(t)}_{2^n}-\mu\circ K\|_{(A_t)_*}
 &=\|(\mu\circ\theta^{(s)}_{2^n}-\mu)\circ K\|_{(A_t)_*}\\
 &\le\|\mu\circ\theta^{(s)}_{2^n}-\mu\|_{(A_s)_*}.
\end{aligned}
\end{equation}
Applying \eqref{eq:hm-state-displacement} to both sides of
\eqref{eq:flow-state-contraction} gives \(t\le s\).
An equivariant Markov map \(A_s\to A_t\) would similarly give
\(s\le t\). Thus, Theorem~\ref{thm:criterion}(2) implies
\[
 L_p(M^{(t)})\simeq L_p(M^{(s)})\quad\Longrightarrow\quad t=s.
\]
The family \(\{L_p(M^{(t)}):0\le t<2\}\) therefore contains
\(2^{\aleph_0}\) pairwise nonisomorphic Banach spaces.

For the upper bound, choose a countable dense sequence in any
separable Banach space. The norms of all
\(\mathbb Q(i)\)-linear combinations of that sequence determine
the normed space on its rational span and hence its completion.
There are countably many such combinations, so there are at most
\(2^{\aleph_0}\) possible norm functions. This bounds the number
of separable Banach isomorphism classes by \(2^{\aleph_0}\).
\end{proof}

Let \(R_\infty\) be the hyperfinite type \(\mathrm{III}_1\) factor.
Its continuous core has center \(\Zc(R_\infty)=\mathbb C\).
For every \(M\) as in Theorem~\ref{thm:criterion}, the map
\[
 \mathbb C\longrightarrow\Zc(M),\qquad
 a\longmapsto a1_{\Zc(M)}
\]
is an equivariant Markov map. A normal unital positive map
\(\mu:\Zc(M)\to\mathbb C\) is equivariant exactly when
\(\mu\circ\theta^M_t=\mu\) for all \(t\in\mathbb R\).
Thus, Theorem~\ref{thm:criterion}(2) gives
\begin{equation}\label{eq:rinfty-class}
 \begin{gathered}
 L_p(M)\simeq L_p(R_\infty)
 \quad\Longleftrightarrow\quad
 (\Zc(M),\theta^M)\text{ admits an invariant normal state},\\
 1\le p<\infty,\qquad p\ne2.
 \end{gathered}
\end{equation}
For \(p=1\) and \(M\) a factor, \eqref{eq:rinfty-class} has the
same flow condition as the cb-isomorphism criterion of
Haagerup and Musat \cite{HM}.
Let \(R_\lambda\) be the hyperfinite type \(\mathrm{III}_\lambda\)
factor, \(0<\lambda<1\). Its periodic flow has an invariant normal
state given by normalized Lebesgue measure on the period interval.
Hence, \eqref{eq:rinfty-class} gives
\(L_p(R_\lambda)\simeq L_p(R_\infty)\).
The complete isomorphisms between these spaces were established by
Haagerup--Rosenthal--Sukochev \cite{HRS}.
Corollary~\ref{cor:continuum-classes} distinguishes the family
\(\{L_p(M^{(t)}):0\le t<2\}\) by the parameter \(t\).
For arbitrary \(M,N\), Theorem~\ref{thm:criterion}(2) gives the
criterion in terms of their two flows of weights.
At \(p=2\), all the spaces \(L_2(M)\) under consideration are
separable infinite-dimensional Hilbert spaces and are therefore
linearly isometric.


\begin{thebibliography}{99}
\bibitem{AH} H. Ando and U. Haagerup, \emph{Ultraproducts of von Neumann algebras},
J. Funct. Anal. \textbf{266} (2014), 6842--6913.
\bibitem{HM} U. Haagerup and M. Musat, \emph{Classification of hyperfinite factors up to completely bounded isomorphism of their preduals},
J. Reine Angew. Math. \textbf{630} (2009), 141--176.
\href{https://doi.org/10.1515/CRELLE.2009.037}{DOI}.
\bibitem{HRS} U. Haagerup, H. P. Rosenthal and F. A. Sukochev,
\emph{Banach embedding properties of non-commutative \(L^p\)-spaces},
Mem. Amer. Math. Soc. \textbf{163} (2003), no. 776.
\href{https://doi.org/10.1090/memo/0776}{DOI}.
\bibitem{HS} U. Haagerup and E. St{\o}rmer,
\emph{Equivalence of normal states on von Neumann algebras and the flow of weights},
Adv. Math. \textbf{83} (1990), 180--262.
\href{https://doi.org/10.1016/0001-8708(90)90078-2}{DOI}.
\bibitem{HuangSukochev2026} J. Huang and F. Sukochev,
\emph{Banach non-isomorphism of type III and semifinite noncommutative
\(L_p\)-spaces}, \href{https://arxiv.org/abs/2610.11055}{arXiv:2610.11055} (2026).
\bibitem{JRX} M. Junge, Z.-J. Ruan and Q. Xu,
\emph{Rigid \(\mathcal{OL}_p\) structures of non-commutative \(L_p\)-spaces
associated with hyperfinite von Neumann algebras}, Math. Scand.
\textbf{96} (2005), 63--95.
\bibitem{ORS} T. Oikhberg, H. P. Rosenthal and E. St{\o}rmer,
\emph{A predual characterization of semi-finite von Neumann algebras},
Contemp. Math. \textbf{335} (2003), 243--245.
\bibitem{Pietsch2007} A. Pietsch, \emph{History of Banach Spaces and Linear Operators},
Birkh\"auser, Boston, 2007.
\bibitem{Raynaud} Y. Raynaud, \emph{On ultrapowers of non commutative \(L_p\) spaces},
J. Operator Theory \textbf{48} (2002), 41--68.
\bibitem{RX} Y. Raynaud and Q. Xu, \emph{On subspaces of non-commutative \(L_p\)-spaces},
J. Funct. Anal. \textbf{203} (2003), 149--196.
\href{https://doi.org/10.1016/S0022-1236(03)00045-4}{DOI}.
\bibitem{Sherman} D. Sherman, \emph{On the structure of isometries between
noncommutative \(L^p\) spaces}, Publ. Res. Inst. Math. Sci. \textbf{42} (2006), 45--82.
\bibitem{Sukochev2000} F. A. Sukochev, \emph{Linear-topological classification
of separable \(L_p\)-spaces associated with von Neumann algebras of type I},
Israel J. Math. \textbf{115} (2000), 137--156.
\href{https://doi.org/10.1007/BF02810584}{DOI}.
\bibitem{Tomatsu} R. Tomatsu, \emph{Ultraproducts of crossed product von Neumann algebras},
Illinois J. Math. \textbf{61} (2017), 275--286.
\end{thebibliography}
\end{document}